# Solving Unit Commitment Problems with Graph Neural Network based Initial Commitment Prediction and Large Neighborhood Search

Linfeng Yang [12], Peilun Li [1] and Jinbao Jian[3]

[1] School of Computer Electronics and Information, Guangxi University, Nanning 530004, China
`ylf@gxu.edu.cn, lipeilun@st.gxu.edu.cn`
[2] Guangxi Key Laboratory of Multimedia Communication and Network Technology, Guangxi University
[3] College of Mathematics and Physics, Center for Applied Mathematics of Guangxi, Guangxi Minzu University, Nanning, Guangxi, 530006, China.

**Abstract.** Unit commitment problem (UCP) is a critical component of power market decision-making. However, its computational complexity necessitates efficient solution methods. In this work we propose a framework to accelerate the solving process of the UCP, and the data collecting process for two distinct graph neural network (GNN) policy. We at first train a *Neural Initial Commitment Prediction* policy to obtain an initial commitment for UCP. Second, a heuristic process is introduced to restore the feasibility of the initial commitment. Third, get the neighborhood based on the initial prediction then neighborhood search to improve the commitment. At last, we train a *Neural neighborhood Prediction* policy to predict the neighborhood of the incumbent commitment at each iteration, continuously optimizing the commitment until the stopping condition is met. This approach produces high-quality initial commitments that can be iteratively refined to meet higher accuracy requirements. The experimental results show that the GNN policies trained on the 80-unit system outperform commercial solvers on a 1080-unit system, and LNS performs better than commercial solver on more complex instances.

## NOMENCLATURE

Indices:
$g$      Index of unit.
$t$      Index of time period.
Constants:
$N$      Total number of units.
$T$      Total number of time periods.
$\alpha_g, \beta_g$      Coefficients of the linear production cost function of unit $g$.
$C_{hot,g}$      Hot startup cost of units $g$.
$C_{cold,g}$      Cold startup cost of units $g$.
$\underline{T}_{on,g}$      Minimum up time of units $g$.
$\underline{T}_{off,g}$      Minimum down time of units $g$.
$\underline{T}_{cold,g}$      Cold startup time of units $g$.
$\overline{P}_g$      Maximum power output of units $g$.
$\underline{P}_g$      Minimum power output of units $g$.
$P_{D,t}$      System load demand in period $t$.
$R_t$      Spinning reserve requirement in period $t$.
$P_{up,g}$      Ramp up limit of unit $g$.
$P_{down,g}$      Ramp down limit of unit $g$.

This work was supported by the National Natural Science Foundation of China (72361003) and the Key Research and Development Program of Guangxi (2023AB01242). (Corresponding author: Linfeng Yang).

$P_{\text{start},g}$   Startup ramp limit of unit $g$.
$P_{\text{shut},g}$   Shutdown ramp limit of unit $g$.
$u_{g,0}$   Initial commitment state of unit $g$.
$T_{g,0}$   Number of periods unit $g$ has been online(+) or offline(-) prior to the first period of the time span (end of period 0).
$U_g$   $max(0, min(T, u_{g,0}(\underline{T}_{\text{on},g} - T_{g,0}))$
$L_g$   $max(0, min(T, u_{g,0}(\underline{T}_{\text{off},g} + T_{g,0}))$

Binary variables:
$u_{g,t}$   Commitment status of unit $g$ in period $t$, equal to 1 if unit $g$ is online in period $t$ and 0 otherwise.
$s_{g,t}$   Startup status of unit $g$ in period $t$, equal to 1 if unit $g$ starts up in period $t$ and 0 otherwise. This variable is applicable solely to the 3-bin UCP.
$d_{g,t}$   Shutdown status of unit $g$ in period $t$, equal to 1 if unit $g$ shuts down in period $t$ and 0 otherwise. This variable is applicable solely to the 3-bin UCP.

Continuous variables
$P_{g,t}$   Power output of unit $g$ in period $t$.
$S_{g,t}$   Startup cost of unit $g$ in period $t$.

# 1  Introduction

The unit commitment problem (UCP) has been receiving significant attentions from both industry and academia. In general, the UCP determine the operation schedule of the generating units at each time period with varying loads under different operating constraints and environments [1]. Numerous methods have been proposed to address the UCP, including Lagrangian relaxation [2], dynamic programming [3], priority list [4] and mixed integer programming (MIP). With advancements in MIP theory and the efficiency of general-purpose branch-and-cut solvers [5], MIP-based approaches are gaining popularity for solving UCPs. Given the NP-hard nature of UCP calculations [6], researchers have focused on improving efficiency through tighter modeling [7-9], and linear programming (LP) relaxation [10]. Recently, machine learning-based approaches have emerged, with several studies exploring their potential for enhancing the UCP solution process [11-14].

## 1.1  LP relaxation for UCP and MIP

Given the computational burden of UCPs, numerous studies have focused on efficiently obtaining near-optimal solutions and accelerating the branch-and-bound (B&B) process. For instance, LP relaxation solutions are leveraged to generate high-quality initial solutions and LNS neighborhoods.

To expedite UCP calculations, Reference [15-17] reconstruct a tight UCP formulation and use heuristics to repair the continuous relaxation variables to obtain high-quality sub-optimal solutions, the neighborhood obtained by the heuristic is then used to iteratively improve the solution. Reference [10] utilize the commitment variables and generation variables in the LP relaxation solution, and fixed the values of some of the variables through a simple heuristic method to reduce the problem size, thus improve the solution efficiency. Reference [18] propose a relaxation-based neighborhood search for LP relaxation and improved relaxation inducement algorithm. Reference [19] select variables in neighborhood by solving the LP relaxation of Local Branch (LB) iteratively and neighborhood search to improve solution for MIP.

## 1.2  ML for solution predicting

Most systems exhibit characteristics, such as the parameters of the generating units, that remain almost completely constant between instances. Machine learning can effectively extract information from previously solved instances and leverage this information to significantly accelerate the solution of similar instances in the future. ML have demonstrated impressive performance in predicting commitments for the UCP. Literature suggests that, for large-scale power systems, learning the solution to the UCP represents a low-hanging fruit that can be effortlessly picked [20].

Reference [21] use support vector machines and *k*-nearest neighbor to predict active transmission constraints can be safely omitted and partial solution as a warm start. Reference [22] propose an E-Sequence to Sequence-based data-driven policy for multiple-sequence mapping samples. Reference [23] proposed a ML based classifier to predict a commitment and restore the feasibility of the predicted commitment.

Due to the permutation invariance that exists in graph neural networks (GNN) representation of MIP [20][24] and GNN based policy has more generalization in the UCP [25], which has spawned a large number of methods that use GNNs to assist MIP in recent years.

Reference [26] encodes the grid topology into a graph and uses a GNN to predict to obtain both active transmission constraints and commitment, then restore the commitment to feasible and reduce problem size base on active transmission constraints. Reference [27] uses graph neural networks to predict the marginal probabilities of each binary variable of the MIP, and then search for the best feasible solution within a properly defined ball around the predicted solution. But search in a properly defined ball is inefficient in ensure the probability and optimality.

### 1.3    ML for solution improving

After predicting the commitment, using a commercial solver with the commitment as an initial solution is a good choice for further improving its quality. While this does not address the time complexity issue of the B&B algorithm [28,29], which remains the cornerstone of commercial solvers. Thus, multiple scholars have proposed methods to accelerate branching, such as imitating Full Strong Branch (FSB) and Local Branch to aid B&B framework and others try to replace B&B with Large Neighborhood Search (LNS). Reference [29] uses GNN networks to learn FSB branching policies to guide the solver to choose the optimal branch in B&B. But, due to its exhaustive search nature, it is still hard for B&B to be useful in large-scale MIPs [28].

LNS has been shown to find high-quality solutions significantly faster than B&B for large-scale MIPs [30-32]. Reference [31] use GNN networks to imitate learning LB policy for the current solution to choose neighborhood and iteratively improving the solution. Reference [32] uses contrastive learning to allow GNN to learn LB strategy, which outperforms ordinary imitation learning as well as LP relaxation-based methods in the tests.

### 1.4    Contributions

Given the existing MIP formulation of the UCP, our focus is on encoding the MIP formulation into a graph structure, rather than relying on the power system topology as in [33]. This approach enables us to exploit the strengths of machine learning techniques tailored for MIPs.

While neural networks offer potential for accelerating UCP, their inability to guarantee feasible commitments and the limitations of current feasibility restoration methods necessitate novel approaches. To address these challenges, we introduce the IP-LNS (Initial Predicting and Large Neighborhood Search) framework, designed to expedite the optimization process while simultaneously generating dataset for training two distinct GNN policies efficiently.

The framework is shown in **Fig.1,** which consists of the following components: (1) a GNN-based policy that predicts generator commitment for the UCP; (2) a restoration process that restore the predicted values to ensure feasible commitments and identifies the less confident parts of the predicted values; (3) get the neighborhood based on the initial prediction then neighborhood search to improve the commitment; and (4) another GNN-based policy that predicts the neighborhood used by LNS, improving the incumbent commitment step by step. The contributions of this paper are summarized as follows:

1) A framework dedicated to UCP encompasses the initial commitment predicting with feasibility restoration and iteratively LNS with adoptive neighborhood size.
2) A tripartite relational graph convolutional neural network(R-GCN) designed to capture the subtle interactions among different constraints, variables, and the objective function.
3) A highly efficient method for collecting datasets that can simultaneously serve *Neural Initial Commitment Prediction* policy and *Neural neighborhood Prediction* policy.

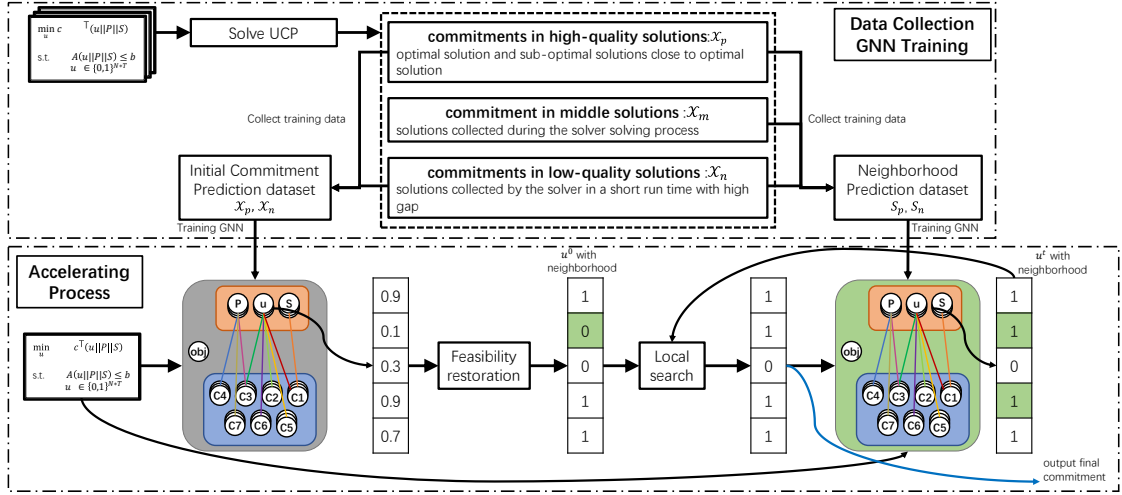

**Fig. 1.** An overview of data collection and our approach is as follows: Initially, we collect high-quality, middle, and low-quality commitments to train both the Neural Initial Commitment Prediction policy and Neural neighborhood Prediction policy. During testing, the input is a UCP formulated as a MIP, and the Neural Initial Commitment Prediction policy outputs the fractional value (between 0 and 1) of each commitment. A heuristic algorithm then restores a feasible commitment and generates a neighborhood (indicated by green). The sub-MIP, which contains only the commitments within the neighborhood (with other commitments fixed to feasible values), is solved using a state-of-the-art MIP solver. This process iterates with the Neural neighborhood Prediction policy predicting the neighborhood based on the incumbent commitment until a runtime limit is reached or other stopping condition is met.

The rest of this paper is organized as follows. Section 2 presents the UCP formulations considered in this work. Section 3 describes the transformation of a UCP into a tripartite graph. While Section 4 models the ML architecture and the proposed methods while the results are discussed in Section 5. Finally, Section VI concludes the paper. Our code and other resources are available at https://github.com/peirisi/IP-LNS.

## 2    Preliminaries

### 2.1   UCP formulations

This paper focuses on classical 1-bin and 3-bin UCP formulations [7]. The 1-bin formulation is simpler, while the 3-bin formulation imposes tighter constraints. If the constraints applied in the 1-bin and 3-bin formulas differ, they are correspondingly labeled. Although this study concentrates on 1-bin and 3-bin models, the proposed approach is adaptable to a broader range of UCP formulations.

**Objective Function:**

The objective function represents the total cost of the power system, which is shown as:

$$min \sum_{g=1}^{N} \sum_{t=1}^{T} [\alpha_g u_{g,t} + \beta_g P_{g,t} + S_{g,t}] \tag{1}$$

where $F_C$ is the total operation cost, the production cost is $\alpha_g u_{g,t} + \beta_g P_{g,t}$, and the startup cost is $S_{g,t}$.

**Startup cost:**

The startup cost $S_{g,t}$ can be expressed as a MIP formulation in 1-bin formulation:

$$\begin{cases} S_{g,t} \geq K_{g,\tau}[u_{g,t} - \sum_{j=1}^{\tau} u_{g,t-j}] \\ S_{g,t} \geq 0, i = 1, \dots, N; t = 1, \dots, T; \tau = 1, \dots, N_{D,g} \end{cases} \tag{2}$$

Where $N_{D,g}$ is a given parameter and $K_{g,\tau}$ models the startup cost as a stepwise function that becomes more accurate as the number of intervals increases:

$$K_{g,\tau} = \begin{cases} C_{hot,g} : \tau = 1, \dots, \underline{T}_{\text{off},g} + T_{cold,g} \\ C_{cold,g} : \tau = \underline{T}_{\text{off},g} + T_{cold,g} + 1, \dots, N_{D,g} \end{cases} \tag{3}$$

With $s_{g,t}$ and $d_{g,t}$ the startup cost $S_{g,t}$ can be expressed as a MIP formulation in 3-bin formulation with:

$$\begin{cases} S_{g,t} \geq C_{\text{hot},g} s_{g,t} \\ S_{g,t} \geq C_{\text{cold},g} \left[ s_{g,t} - \sum_{\tau=\max(t-T_{\text{off},g}-T_{\text{cold},g},1)}^{t-1} d_{i,\tau} - f_{\text{init},g,t} \right] \end{cases} \quad (4)$$

Where $f_{\text{init},g,t} = 1$ when $t - \underline{T}_{\text{off},g} - T_{\text{cold},g} \leq 0$ and $max(-T_{g,0}, 0) < |t - \underline{T}_{\text{off},g} - T_{\text{cold},g} - 1| + 1$, elsewise $f_{\text{init},g,t} = 0$.

**Unit generation limits:**

$$u_{g,t} \underline{P}_g \leq P_{g,t} \quad (5)$$

$$P_{g,t} \leq u_{g,t} \bar{P}_g \quad (6)$$

**Power balance constraint:**

$$\sum_{g=1}^{N} P_{g,t} - P_{D,t} = 0 \quad (7)$$

**System spinning reserve requirement:**

$$\sum_{g=1}^{N} u_{g,t} \bar{P}_g \geq P_{D,t} + R_t \quad (8)$$

**Ramp rate limits:**

Ramp rate limits constraints of the 1-bin formulation as follows:

$$P_{g,t} - P_{g,t-1} \leq u_{g,t-1} P_{up,g} + (u_{g,t} - u_{g,t-1}) \quad (9)$$

$$P_{g,t-1} - P_{g,t} \leq u_{g,t} P_{down,g} + (u_{g,t-1} - u_{g,t}) \quad (10)$$

With $s_{g,t}$ and $d_{g,t}$ the Ramp rate limits constraints of the 3-bin formulation as follows:

$$P_{g,t} - P_{g,t-1} \leq u_{g,t}(P_{up,g} + \underline{P}_g) - u_{g,t-1} \underline{P}_g + s_{g,t}(P_{start,g} - P_{up,g} - \underline{P}_g) \quad (11)$$

$$P_{g,t-1} - P_{g,t} \leq u_{g,t-1}(P_{down,g} + \underline{P}_g) - u_{g,t} \underline{P}_g + d_{g,t}(P_{shut,g} - P_{down,g} - \underline{P}_g) \quad (12)$$

**Minimum up/down time constraints:**

Minimum up/down time constraints of the 1-bin formulation as follows:

$$u_{g,t} - u_{g,t-1} \leq u_{g,\tau}, \tau \in [t+1, min\{t + \underline{T}_{on,g} - 1, T\}] \quad (13)$$

$$u_{g,t-1} - u_{g,t} \leq 1 - u_{g,\tau}, \tau \in [t+1, min\{t + \underline{T}_{off,g} - 1, T\}] \quad (14)$$

With $s_{g,t}$ and $d_{g,t}$ the Minimum up/down time constraints of the 3-bin formulation as follows:

$$\sum_{\varpi=\max(t-\underline{T}_{on,g},0)+1}^{t} s_{g,\varpi} \leq u_{g,t}, t \in [U_g + 1, \ldots, T] \quad (15)$$

$$\sum_{\varpi=max(t-\underline{T}_{off,g},0)+1}^{t} d_{g,\varpi} \leq 1 - u_{g,t}, t \in [L_g + 1, \ldots, T] \quad (16)$$

**Initial status of units:**

$$u_{g,t} = u_{g,0}, t \in [1, \ldots, U_g + L_g] \quad (17)$$

**State constrains:**

State constraints were introduced in the 3-bin formulation to link the startup, shutdown, and state variables.

$$s_{g,t} - d_{g,t} = u_{g,t} - u_{g,t-1} \quad (18)$$

## 2.2 Large Neighborhood Search

LNS is a heuristic algorithm that starts with an initial solution and then iteratively picks a neighborhood and performs a neighborhood search to improve the solutions until a runtime limit is exceeded or other stopping condition is met. Compared to B&B, LNS is more effective in improving the objective value, especially on difficult instances [31,34,35]. Compared to other local search methods, LNS explores a large neighborhood in each step and thus, is more effective in avoiding local minima.

Let $\delta^k = \{M, u^k\}, u^k \in \{0,1\}^{N \times T}$ denote a UCP $M$ with commitment $u^k$ from iteration $k$. The initial commitment $u^0$ in this paper is predicted by a *Neural Initial Commitment Prediction* policy. The neighborhood $a^k \in \{0,1\}^{N \times T}$, predicted by a *Neural neighborhood Prediction* policy, represents the neighborhood at iteration $k$. A value of 0 indicates that the corresponding commitment variable will be fixed to the value in $u^k$, constructing a sub-MIP with fewer variables. The sub-MIP is then solved to obtain the commitment solution $u^{k+1}$, leading to the next iteration at $k+1$ with a nonincreasing objective value.

## 3 Tripartite Graph Representation for UCPs

In this section, we use tripartite graphs to represent UCPs. The use of tripartite graphs to represent MIPs was proposed in [36], while [29] introduced a bipartite graph representation of the MIP problem, which complicates targeted optimization due to the generalized nature of MIPs. We propose a tripartite R-GCN to model the complex interactions within the UCP, capturing relationships between constraints, variables, and the objective function.

The transformation from the 1-bin UCP to a tripartite graph is illustrated in **Fig.2**. For simplicity the 1-bin UCP is used and the names of the constraints in **Fig.2** are represented using shorthand:

C1: Startup cost
C2: System spinning reserve requirement
C3: Unit generation limits
C4: Power balance constrains
C5: Initial status of units
C6: Minimum up/down time constraints
C7: Ramp rate limits

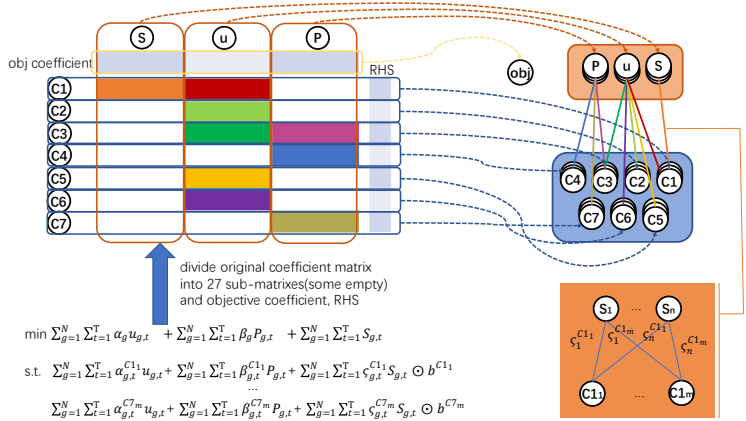

**Fig. 2.** Transformation of a 1-bin UCP into a tripartite graph involves integrating variables and constraints within the MIP modeling framework to facilitate the division of the coefficient matrix in the tripartite graph transformation. Each connection between a type of variable and a type of constraint represents a sub-matrix of the coefficient matrix, distinguished by color coding.

Let $G = (V \cup C \cup O, E^{VC} \cup E^{Vo} \cup E^{Co})$ denote a 1-bin UCP tripartite graph, where $V = V^u \cup V^P \cup V^S$ denotes 3 disjoint sets, each of these represents a type of variable, $V^v = \{V_1^v, V_2^v, \ldots, V_{n_v}^v\}, v \in \{u, P, S\}$, each of which corresponds to a variable; $C = C^{C1} \cup C^{C2} \cup \ldots \cup C^{C7}$, 7 kinds of constraints set in set $C$, $C^c = \{C_1^c, C_2^c, \ldots, C_{m_c}^c\}, c \in \{C1, C2, \ldots, C7\}$, each of which corresponds to a constraint, and objective function vertex $O$ corresponds to the objective function. $vc, vo, co$ represents the set of edges between variables and constraints, variables and objective function, constraints and objective function.

$E^{VC} = \{\mathcal{E}^{vc} | v \in \{u, P, S\}, c \in \{C1, C2, \ldots, C7\}\}$ where there exist $3 \times 7$ kinds of categories corresponding to the edges between the set of 3 variables and the set of 7 constraints, but edges do not exist for every variable with every constraint. An edge $(i, j) \in \mathcal{E}^{vc}, v \in \{u, P, S\}, c \in \{C1, C2, \ldots, C7\}$ exists only if the coefficient of the variable $v_i$ in the constraint $c_j$ is not zero, and the features of the edge $\mathcal{E}_{i,j}^{vc}$ correspond to the coefficients in the coefficient matrix;

$E^{Vo} = \{\mathcal{E}^{uo}, \mathcal{E}^{Po}, \mathcal{E}^{So}\}$ where there exist 3 kinds of edges corresponding to the edges of the set of 3 variables to the objective function. An edge $(i) \in \mathcal{E}^{vo}, v \in \{u, P, S\}$ exists only if variables $v_i$ with non-zero coefficients in the objective function, and the features of the edges $\mathcal{E}_i^{vo}$ correspond to the coefficients of the variables in the objective function;

$E^{Co} = \{\mathcal{E}^{C1o}, \mathcal{E}^{C2o}, \ldots, \mathcal{E}^{C7o}\}$ where there are 7 kinds of edges corresponding to the 7 constraint sets to the objective function. An edge $(j) \in \mathcal{E}^{co}, c \in \{C1, C2, \ldots, C7\}$ exist only if this constraint $c_j$ is active in the solution. The features of the edges $\mathcal{E}_j^{co}$ correspond to the RHS of the constraints. For the initial commitment prediction task, it is the LP relaxation solution and for the neighborhood prediction task it is the incumbent solution.

## 4 Framework

1) A GNN-based architecture predicts each generator commitment variables in fractional form.
2) Heuristics are used to obtain a feasible commitment and neighborhood based on the predicted fractional values.
3) Solving sub-MIP based on the neighborhood;
4) Another GNN-based architecture iteratively predicts neighborhood base on incumbent solution and solver solve the sub-MIP until a runtime limit is exceeded or other stopping condition is met.

### 4.1 Policy Network

We learn a policy $\pi_\theta$ contract by MLP and R-GCN, parameterized by learnable weights $\theta$. The input of policy is a tripartite graph $G = (V \cup C \cup O, E^{VC} \cup E^{Vo} \cup E^{Co})$ that represent a UCP. The policy is a three-layer structure describe as follow:

**The input layer** consisting of a regularization layer, fully connected layers, and activation functions, with parameters not shared between different kinds. In this case, the LP relaxation feature dimensions of the variables need to skip the normalization layer at the beginning, because the normalization layer will disable the LP relaxation solution. The corresponding experiments will be presented in Section 5.

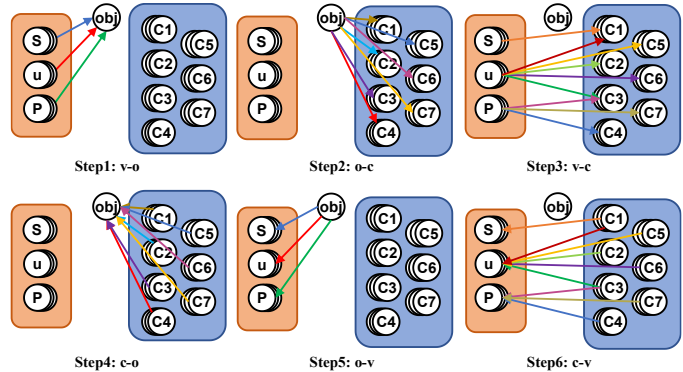

**Fig. 3.** Information transition flow in the trigraph convolutional layer

**GCN layer**, the information of neighboring nodes is transitioned after multilayer GCN. The information transition flow is shown in **Fig. 3**, and take the form:

Variables to objective:
$$O \leftarrow f_o(O, ln(\sum_{v \in \{u,P,S\}} \sum_{(i) \in \mathcal{E}^{vo}} f_{vo}(V^v{}_i, O, \mathcal{E}_i^{vo})))$$

Objective to Constrain:
$$C_j^c \leftarrow f_c(C_j^c, ln(f_{oc}(O, C_j^c, \mathcal{E}_j^{co})), c \in \{C1, C2, \ldots, C7\}$$

Variables to Constrain:
$$C_j^c \leftarrow f_c(C_j^c, ln(\sum_{v \in u,P,S} \sum_i^{(i,j) \in \mathcal{E}^{vc}} f_{vc}(V_i^v, C_j^c, \mathcal{E}_{i,j}^{vc}))), c \in \{C1, C2, \ldots, C7\}$$

Constrains to Objective:
$$O \leftarrow f_o(O, ln(\sum_{c \in \{C1,C2,\ldots,C7\}} \sum_{(j) \in \mathcal{E}^{co}} f_{co}(C_j^c, O, \mathcal{E}_j^{co})))$$

Objective to Variable:
$$V_i^v \leftarrow f_v(V_i^v, ln(f_{ov}(O, V_i^v, \mathcal{E}_i^{vo}))), v \in \{u, P, S\}$$

Constrains to Variable:
$$V_i^v \leftarrow f_v(V_i^v, ln(\sum_{c \in \{C1,C2,\dots,C7\}} \sum_j^{(i,j) \in \mathcal{E}^{vc}} f_{cv}(C_j^c, V_i^v, \mathcal{E}_{i,j}^{vc}))), v \in \{u, P, S\}$$

For all $v \in \{u, P, S\}, c \in \{C1, C2, \dots, C7\}$ where the functions $f$ are MLP with ReLU activation functions. Different subscripts are used to distinguish various kinds, and $ln$ represents the layer normalization function.

**Output layer** consists of fully-connected layers and an activation function that reduces the dimensional output of the commitment variables. Other variables are discarded.

It is worth noting that both predict initial commitment policy and neighborhood policy in this paper share the same GNN policy and output dimensionality. Only the dimensionality of the input data is different.

### 4.2 Predicting Initial Commitment

In this section, our goal is to train a *Neural Initial Commitment Prediction* policy using supervised learning to predict the generator commitments. Given the computational expense of solving the UCP, we employ a contrastive learning approach to train the policy efficiently. A set of high-quality commitments are collected as positive samples, while a set of lower-quality commitments serve as negative samples. By contrasting these samples, the policy learns to distinguish between effective and suboptimal commitments. We have adapted the input features proposed in [27] better align with *Neural Initial Commitment Prediction* policy.

### 4.3 Feasibility Restoration

While graph neural networks can predict high-quality commitments, directly fixing certain commitment variables based on these predictions may introduce the risk of infeasibility. The method of using a feasibility pump[26,23] to restore an infeasible commitment to a feasible commitment is usually effective when the predicted commitment is of high quality, but the time taken to repair a feasible commitment is difficult to guarantee.

Thus, we use a feasibility repair method based on heuristics to ensure the feasibility of the commitment without consuming too much time. And stability metrics are to be preserved to avoid meaningless neighborhoods. Feasibility restoration methodologies improved from [16]:

Due to the complexity of the UCP, there is always the possibility that there is no feasible commitment for feasibility restoration, such as turning on too many units resulting in a lower bound of generation turned on that exceeds the required load. However, it is worth noting that this scenario occurs rarely in experimental settings.

Therefore, at the end it is necessary to check whether the commitment is feasible or not, and if it is not feasible then the feasibility pump method is used to find a feasible commitment closest to the current infeasible commitment. The feasibility pump can find a feasible commitment based on current infeasible commitment $u^*$, can be formulated as follow:

$$\min \sum_{i=1}^{N} \sum_{t=1}^{T} |u_{i,t} - u_{i,t}^*|$$

$$s.t. \begin{cases} (2-3)(5-8)(9-10)(13-14)(17), 1-\text{bin} \\ (4)(5-8)(11-12)(15-16)(17)(18), 3-\text{bin} \end{cases}$$

Solving this MIP with Gurobi solver, a feasible commitment is guaranteed, albeit with additional time. It is worth noting that the aforementioned situation leading to infeasibility is rarely triggered. Gurobi is configured to exit upon finding the first feasible commitment, thereby quickly restoring the feasibility of commitment.

### 4.4 Local Search

The commitments variables are fixed to the values taken after restoration, and a neighborhood can be calculated by the processed fractional prediction. Algorithm 1 presents the details of local search process. For high-confidence predictions, we set constraints to fix their values; for low-confidence predictions, we set the solver's initial values to local search improve the commitment.

**Algorithm 1** Local Search Algorithm
**Parameter:** lower/upper threshold $\{lt, ut\}$
**Input:** predict commitment after extra process in restore $u \in (0,1)^{(N \times T)}$, restored commitment $u^* \in \{0,1\}^{(N \times T)}$, UCP $M$, commitment variables x

1: **for** $g = 1 : N, t = 1 : T$ **do**
2:   **if** $u_{g,t} <= lt$ or $u_{g,t} >= ut$ **then**
3:     $M$ add constraint $(x_{g,t}, u_{g,t}^*)$
4:   **else**
5:     $M$ set start $(x_{g,t}, u_{g,t})$
6:   **end if**
7: **end for**
6: solve $M$

When restoring the minimum on/off time or ramp constraints, there will be some commitments that are set to on due to other commitments which usually has low reliability so it is easy to be chosen as a neighbor while exists meaningless neighborhoods.

As shown in **Fig. 4**, the yellow variables must remain in the 'on' state. These variables are included in the neighborhood but their values cannot change; if any of these values were to shift from 1 to 0, the commitment would become infeasible. This results in meaningless neighborhoods. Thus, a sequence of consecutive on-state predictions is set to the highest value among them to avoid meaningless neighborhoods.

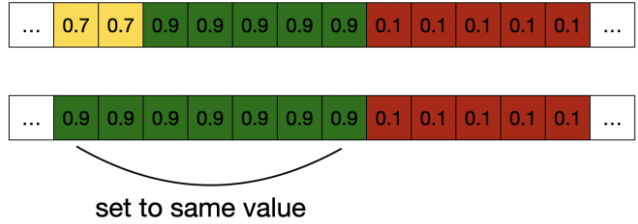

**Fig. 4.** A part of unit $g$ schedule, the minimum $\underline{T}_{on,g}$ is 7. Green indicates prediction values in this period is higher that $ut$, so unit $g$ fixed to 'on' in this period, red indicates prediction values in this period is lower than $lt$, so unit $g$ fixed to 'off' in this period, and yellow indicates unit $g$ that were restored to 'on' during the feasibility restoration process but prediction values between $lt$ and $ut$.

### A. LNS

While the *Neural Neighborhood Prediction* policy's features were expanded upon, additional features were incorporated from [26,31]. After obtained a neighborhood based on the prediction, there is no more information for another neighborhood. Thus, we have a second policy to predict neighborhood base on incumbent commitment. We use contrastive learning to train GNN networks thus we need positive neighborhoods and negative neighborhoods[32], however, unlike it, we did not choose to imitate LB, because for a high-quality prediction, even a small neighborhood can find the optimal solution through the LB algorithm, and it is difficult to generate multiple pieces of data based on the initial commitment for a given UCP, which may result in the GNN over-learning the information of the UCP instead of focusing on how to improve any given commitment.

Thus, we collect high-quality commitments $\mathcal{X}_p$, intermediate commitments $\mathcal{X}_m$ and low-quality commitments $\mathcal{X}_n$. For any intermediate commitment $\dot{x} \in \mathcal{X}_m$, we will generate positive neighborhoods $S_p$ and negative neighborhoods $S_n$ for it. Its positive neighborhoods $S_p$ will retain only those differences between the high-quality commitments and the intermediate commitment that do not exceed $\alpha_{size} \times |\dot{x}|$ and the resulting improvement of neighborhood search is greater than $\alpha_p$ times the best neighborhood improvement in positive neighborhoods, $\alpha_{size}$ is a ratio parameter that exclude oversized neighborhoods, where $sum(\star)$ denotes the summation of $\star$, xor denotes the exclusive OR operation where differing value result in 1 and identical value result in 0, $|\star|$ denotes the length of $\star$. The number of negative neighborhoods is $\alpha_c |S_p|$ where $\alpha_c$ is the ratio parameter of positive and negative sample size. Take $x_n$ out from $\mathcal{X}_n$ iteratively, the neighborhood is $xor(\dot{x}, x_n)$, if the resulting improvement of neighborhood search is less than $\alpha_n$ times the best positive neighborhood search improvement, it will be seen as a negative neighborhood. The iteration continues until a sufficient number of negative neighborhoods have been collected; otherwise, we will perturb neighborhood by 5% and this percentage will increase by 5% per iteration until it reaches 100%, after which only the increase will stop but the iteration will continue. $\alpha_p, \alpha_{size}, \alpha_c, \alpha_n$ are set to 0.6, 0.2, 10 and 0.05, respectively, in experiments.

---
**Algorithm 2 LNS**

Parameter: UCP $M$, $u^0$ commitment before LNS, $u^k$ commitment in step $k$, $a^k$ neighborhood in step $k$, $\pi_\theta$ Neural neighborhood Prediction policy, $\zeta$ neighborhood size ratio.

1 : $gd = 1^{(N \times T)}$, $ld = 1^{(N \times T)}$, $\zeta = 0.2$
2 : for $k = 1 : MAX\_step$ do
3 :    $p^k := \pi_\theta(\delta^k)$
4 :    $p^k = p^k \times gd \times ld$
5 :    $a^k =$ greedy select$(p^k, \zeta) \cup$ row neighborhood$(u^{k-1})$
6 :    $u^k := M$ neighborhood search$(u^{k-1}, a^k)$
7 :    $\zeta =$ adaptive size $(u^{k-1}, u^k)$
8 :    $gd, ld =$ weight descend$(u^{k-1}, u^k, a^k)$
9 :    if stopping condition met then
10:       return $(u^k)$
11:    end if
12: end for
13: return $(u^{MAX\_step})$

---

**Adaptive neighborhood size:** The size in the neighborhood is $\zeta|N \times T|$ determined by a ratio factor $\zeta$, if the objective value solved by neighborhood in this iteration has no improvement, we must increase the neighborhood to increase the search range, then $\zeta = min\{\zeta \times \psi_l, \varphi_l\}$ $\psi_l > 1$ and $\varphi_l$ is a constant to prevent the neighborhood too large to solve. Vice versa if the objective value has improvement, we can shrink the neighborhood reduces the extra burden due to previous increasing, then $\zeta = max\{\zeta \times \psi_u, \varphi_u\}$ $\psi_u < 1$ and $\varphi_u$ is a constant to prevent the neighborhood too small to be ineffective. $\psi_l, \psi_u, \varphi_l, \varphi_u$ are set to 1.1, 0.8, 0.3, 0.1 respectively, in experiments.

**Row neighborhoods:** Usually the first and last few periods of a segment of a generator with consecutive turn-on are a good neighborhood and have been proved in[16]. After each round of neighborhood selection, the incumbent's row neighborhood is added as the neighborhood used for the subsequent neighborhood search.

**Sampling actions:** We greedily select those with the highest prediction value as neighborhood, and like [32], it may select the same neighborhood due to there is only minor difference between commitment to commitment. The entire process is outlined in Algorithm 2. Where the stopping condition is typically triggered by prolonged failure to improve the commitment or by reaching the runtime limit.

So, we use a weight descend method that reduces the weight of commitment variables that have been predicted to be neighbors in the greedy selection. Let $(p_{0,0}, \cdots, p_{N,T}) := \pi_\theta(\delta^k)$ be the commitment fractional prediction output by the policy. Each fractional prediction $p_{g,t}$ is multiplied by two weights before being greedily selected to reduce its weight: the global descending parameter $gd_{g,t}$ and last round descending parameter $ld_{g,t}$. The initial value both are 1, if a prediction $p_{g,t}$ is selected in the neighborhood, then $gd_{g,t} = gd_{g,t} \times \psi_{gd}$ and $ld_{g,t} = \psi_{ld}$, if the commitment variable changed in next step, there will be more descending $gd_{g,t} = gd_{g,t} \times \varphi_{gd}$ and $ld_{g,t} = ld_{g,t} \times \varphi_{ld}$. The global descending parameter reducing the weight in every subsequent round, while the last round descending parameter only reduce the following round. This ensure that the neighborhood will not repeat same. $\psi_{gd}, \psi_{ld}, \varphi_{gd}, \varphi_{ld}$ are set to 0.9, 0.5, 0.8, 0.01 respectively, in experiments.

## 5 Numerical Experiments

### 5.1 Experimental Settings

**Generation system setting:** The generators are not specifically designed to match a particular power system. Thus, thermal generators within the system share identical parameters, meaning each type of generator appears multiple times. The base 8-unit data are taken from [7]. The 60-, 80-, 132-, and 1080-unit test systems are created by duplicating the base 8-unit data.

**Demand setting:** There are 365 loads from the first day of the year to the 365th day of the year come from https://github.com/GridMod/RTS-GMLC. In order to match the loads to the systems used, the loads for the 365 days are mapped to the capacity of the corresponding systems. Each day if the end is 3 it is used as a validation set, if the end is 7 it is used as a test set and the rest is used as a training set. And each load will be perturbed to expand the dataset.

**Device setting:** The data was collected using 3-bin modeling and solved with Gurobi 11.0. The training set was gathered on a system running Ubuntu Server 22.04, equipped with dual Intel Xeon Gold 6140 CPUs @ 2.30GHz, an NVIDIA GeForce RTX 2080 Ti, and 64 GB of memory. Because the server's CPU has a low clock speed, the computation time is excessively long. We conducted the tests on a Windows 11 PC with an Intel Core i7-8700 CPU running at 3.20GHz and 32GB of memory.

**Data collection setting:** For each UCP, the data was solved using 3-bin formulations. When generating data, the model was adapted to either 1-bin or 3-bin based on the requirements. For high-quality commitments set, the MIPGap was set to 1e-7, PoolSearchMode to 2, PoolSolutions to 10, and TLE to 60 seconds. During the solving process, commitments were collected through a callback function as middle commitments set, retaining only those with a relative gap greater than 1e-5 from the average of high-quality commitments set. For low-quality commitments set, the MIPGap was set to 1e-2, PoolSearchMode to 2, PoolSolutions to 90, and TLE to 4 seconds. Data generation was performed in parallel using 50 threads to produce 50 datasets.

All policies are obtained based on 80-unit system training. For initial commitment prediction task, the 80-unit system had 5828 examples as a training set, 738 examples as a validation set. For neighborhood prediction task the 80-units system had 10200 examples as a training set, 1370 examples as a validation set;

**Benchmark setting:** The latest Gurobi solver 11.0 was used with the following settings: Threads set to 1. All other settings were kept at their default values.

**Time calculate setting:** For the UCP, the system size is typically fixed, so only the input loads change on the same system. Based on this, Gurobi can build the model in advantage and then adjust the coefficient matrix as well as the RHS according to the daily loads, which saves a lot of modeling time, so the Gurobi solution time is not calculated for the problem modeling time. It is similar for tripartite graph modeling, where only some of the elements in the matrix need to be adjusted, and therefore tripartite graphs are also not calculated as time spent modeling the graph.

## 5.2 Initial commitment Prediction

The quality of LP relaxation commitments obtained from UCP formulations with different tightness varies, the tighter the formulation has the higher quality relaxation solutions that give better features to the GNN. However, in addition to the effect of LP relaxation, the tightness of the UCP formulation also affects the GNN policies.

**Table 1.** Initial objective value gap (**denote in percentage %**) of different tightness UCP formulation using LP relaxations as features and not using LP relaxations as features on systems of different sizes. The gap is the difference between the primal bound $v$ and the objective value obtained with Gurobi TLE set to 60s $v^*$, defines as $\frac{v-v^*}{v^*}$. A negative value indicates that the predicted commitments have a lower average objective function value than the one obtained by Gurobi within 60 seconds.

| method | LP | 60-unit | 80-unit | 132-unit | 1080-unit |
|---|---|---|---|---|---|
| 1-bin | No | 5.978 | 5.181 | 8.049 | 8.531 |
|  | Yes | 4.237 | 4.195 | 3.193 | 2.283 |
| 3-bin(-) | No | 3.171 | 3.404 | 1.516 | 1.684 |
|  | Yes | 1.307 | 1.289 | 1.132 | 0.048 |
| 3-bin | No | 1.497 | **1.128** | 2.702 | 0.884 |
|  | Yes | **1.480** | 1.337 | **0.912** | **-0.414** |

In **Table 1**, "1-bin" or "3-bin" denotes that the MIP form of UCP is modeled by the 1-bin or 3-bin UCP formulation; "3-bin(-)" indicates the MIP form of UCP that replaces the tighter ramp rate limit constraints of the 3-bin UCP formulation with the looser ramp rate limit constraints from 1-bin UCP formulation. "Yes" or "No" indicates whether LP relaxation is used as a feature. All policies are trained only on 80-unit system but tested on all 4 systems.

According to the results shown in **Table 1**, tighter UCP formulations tend to have lower gap. For test systems similar to the training system, such as 60-unit system, 80-unit system, the GNN policy alone produces fairly good outputs, with little to no difference observed in some cases when using LP fea-tures compared to not using them. However, for systems with significant differences, such as the 132-unit system,1080-unit system, LP relaxation as a feature significantly enhances the predicted initial commitment, 3-bin UCP formulation with LP feature even have lower objective function value than Gurobi within 60s.

### 5.3 Restore

We aim to evaluate the performance of our proposed feasibility restoration method against the feasibility pump on both high-quality and low-quality predictions.

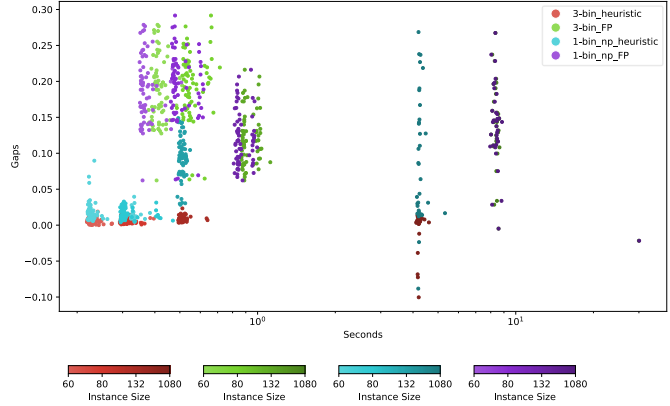

**Fig. 5.** Gap between the 60s Gurobi objective value and the objective value after restoration by different feasible restoration methods on different sized systems as well as on different quality predictions. 72 instances for 60-unit,80-unit,132-unit system and 36 instances for 1080-unit system.

In **Fig. 5**, "3-bin" indicates that the commitment to be restored is predicted by the policy trained on 3-bin and represents high quality predictions; "1-bin_np" indicates that the commitment to be restored is predicted by the policy trained on 1-bin without LP feature representing low quality predictions; "heuristic" indicates that the feasible restoration method used is the previously proposed and "FP" indicates the use of feasibility pump methods. For high-quality predictions, our proposed feasibility restoration method not only requires less time but also yields initial commitments with lower objective function values. The same conclusion applies to low-quality predictions, where the restored commitment after our proposed feasibility restoration method achieves lower objective function values compared to the commitments obtained by the feasibility pump on high-quality predictions. The feasibility restoration method proposed earlier is more competitive than the feasibility pump on all examples regardless of the quality of the predictions.

### 5.4 Ablation study

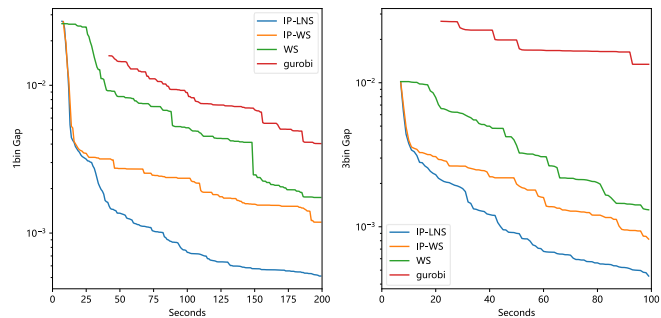

**Fig. 6.** Average gap (relative distance to the optimum among all methods) of the different methods on a 1080-unit system, left 1-bin example, right 3-bin example.

We tested the effectiveness of the individual components in IP-LNS compared to the state-of-the-art commercial solver. "IP-LNS" indicate that predicting the initial commitment by *Neural Initial Commitment Prediction* policy followed by a local search based on the predicted values to determine the neighborhood, then iteratively using the *Neural neighborhood Prediction* policy to predict the neighborhood for Large Neighbor Search; "IP-WS" indicate that predicting the initial commitment by *Neural Initial Commitment Prediction* policy followed by a local search based on the predicted values to determine the neighborhood, the incumbent commitment is then set as the start to warm start Gurobi to further improve the commitment; "WS" indicate that let the predicted initial commitment as a start to warm start Gurobi solver; "Gurobi" indicate that use only Gurobi to solve the MIP.

As shown in **Fig. 6**, each component of IP-LNS plays a key role. In the 3-bin UCP formulation, the average gap of the initial commitment predicted by *Neural Initial Commitment Prediction* policy is even lower than the commitment obtained by Gurobi within 100s. The LNS method has better performance compared to Gurobi's B&B framework in both 1-bin and 3-bin UCP formulations, especially in the 1-bin UCP formulation where the solution processing is less efficient. Neighborhood search based on the initial prediction plus LNS has the best performance in the early stage of the solution processing.

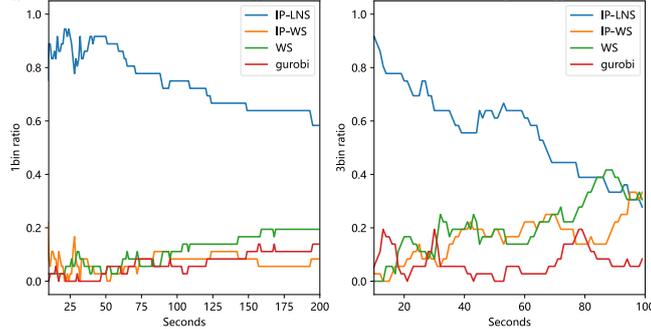

**Fig. 7.** best performance rate (The ratio of instance with objective values that are optimal among all methods, the higher the better) 1-bin left, 3-bin right. The first few seconds yielded no commitments, leading to the truncation of results.

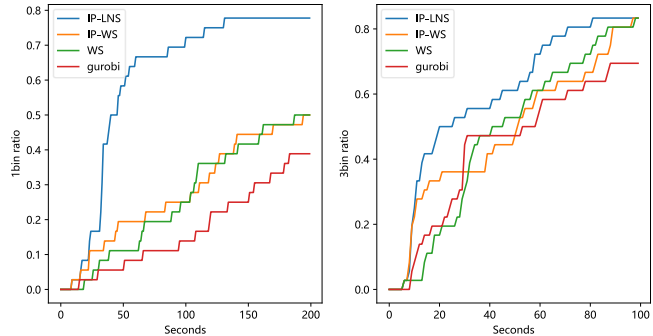

**Fig. 8.** survival rate (the ratio of instance with gap below 1e-3, the higher the better) 1-bin left, 3-bin right.

**Fig. 7** shows the 1-bin UCP formulation and 3-bin UCP formulation best performance rate, LNS has better performance than Gurobi in the early stage. However, as the solution time increases, Gurobi's B&B framework gradually catches up, especially in the more efficient 3-bin UCP formulation. This is because LNS cannot guarantee optimal solutions, while the B&B framework can ensure global optimality through branch and bound.

**Fig. 8** shows the 1-bin and 3-bin UCP formulation primal gap to reach a survival rate of 1e-3, LNS has a better performance than Gurobi in the early stage. LNS clearly has a greater advantage over Gurobi on relatively less tight 1-bin UCP formulation.

# 6 VI. Conclusion

This paper introduces the IP-LNS framework to expedite UCP solving. The framework integrates GNN-based prediction, feasibility restoration, local search, and iterative neighborhood learning to achieve significant speedup without compromising solution quality. Experimental results demonstrate superior performance compared to commercial solvers across different power system scales. While the IP-LNS framework effectively accelerates UCP, its combination with Branch-and-Bound algorithms is envisioned to further enhance its applicability and solution optimality.